# Level sets of the stochastic wave equation driven by a symmetric Lévy noise

DAVAR KHOSHNEVISAN[1] and EULALIA NUALART[2]

[1]*Department of Mathematics, University of Utah, 155 S. 1400 E. Salt Lake City, UT 84112-0090, USA. E-mail: davar@math.utah.edu*
[2]*Institut Galilée, Université Paris 13, 93430 Villetaneuse, France. E-mail: nualart@math.univ-paris13.fr*

We consider the solution $\{u(t,x); t \geq 0, x \in \mathbf{R}\}$ of a system of $d$ linear stochastic wave equations driven by a $d$-dimensional symmetric space-time Lévy noise. We provide a necessary and sufficient condition on the characteristic exponent of the Lévy noise, which describes exactly when the zero set of $u$ is non-void. We also compute the Hausdorff dimension of that zero set when it is non-empty. These results will follow from more general potential-theoretic theorems on the level sets of Lévy sheets.

*Keywords:* level sets; Lévy noise; potential theory; stochastic wave equation

## 1. Introduction and main results

Consider the solution $u(t,x) := (u_1(t,x), \ldots, u_d(t,x))$ to the following system of $d$ linear stochastic wave equations:

$$\frac{\partial^2 u_i}{\partial t^2}(t,x) = \frac{\partial^2 u_i}{\partial x^2}(t,x) + \dot{L}_i(t,x), \qquad t \geq 0, x \in \mathbf{R}, \tag{1.1}$$

with initial condition $u_i(0,x) = \partial_t u_i(0,x) = 0$, where $i$ ranges in $\{1,\ldots,d\}$, and $\dot{L} := (\dot{L}_1, \ldots, \dot{L}_d)$ is a (totally scattered) $d$-dimensional Lévy noise on $\mathbf{R}^2$ with Lévy exponent $\Psi$.

Stochastic PDEs (SPDEs) that are driven by (non-Gaussian) Lévy noises are beginning to receive some attention in the literature. We mention, for example, the work of Mueller [28], who investigates nonlinear heat equations in high space dimensions that are driven by multiplicative space-time stable Lévy noises; see also Mueller, Mytnik and Stan [29] who consider heat equations with a multiplicative space-dependent stable Lévy noise. The Cauchy problem in high space dimensions for the linear wave equation driven by a space-time Lévy noise is treated in Øksendal, Proske and Signahl [30].

We now return to the SPDE (1.1) and define the *light cone* at $(t,x) \in \mathbf{R}_+^2$ to be the set

$$\mathscr{C}(t,x) := \{(s,y) \in \mathbf{R}_+^2 : 0 \le s \le t \text{ and } x-(t-s) \le y \le x+(t-s)\}. \tag{1.2}$$

Then, $\frac{1}{2}\mathbf{1}_{\mathscr{C}(t,x)}$ is the Green function for the wave operator. We follow the approach developed by Walsh [35] for space-time white noise and can write the solution of (1.1), in mild form, as

$$u(t,x) = \frac{1}{2}\int_{\mathbf{R}_+^2} \mathbf{1}_{\mathscr{C}(t,x)}(s,y) L(\mathrm{d}s\,\mathrm{d}y) \qquad \text{for all } t \ge 0, x \in \mathbf{R}, \tag{1.3}$$

where $L$ denotes the Lévy sheet that corresponds to $\dot{L} := (\dot{L}_1, \ldots, \dot{L}_d)$ (see Section 2.1), and the stochastic integral is a Wiener-type integral with respect to $L$ (see Section 2.3).

The aim of this paper is to study the geometry of the zero set

$$u^{-1}\{0\} := \{(t,x) \in \mathbf{R}_+ \times \mathbf{R} : u(t,x) = 0\} \tag{1.4}$$

of the solution to (1.1). In particular, we seek to know when $u^{-1}\{0\}$ is non-void. In order to describe when $u^{-1}\{0\} \neq \varnothing$, we first analyze the zero set of a Lévy sheet. This will be done in a more general setting. Then, we devise a 'comparison principle' to relate the solution of the SPDE (1.1) to theorems about additive Lévy processes developed earlier in [22, 24, 25]. This comparison method might appear somewhat roundabout, but, as it turns out, it is not so easy to work directly with the solution to (1.1).

In order to state the main result of this paper, we introduce the following regularity conditions on the characteristic exponent $\Psi$ of the noise $\dot{L}$; see (2.1) below for a precise definition of $\Psi$.

A1. $\dot{L}$ is *symmetric*. That is, $\Psi(\xi)$ is real and non-negative for all $\xi \in \mathbf{R}^d$.
A2. $\Phi(\lambda) < \infty$ for all $\lambda > 0$, where $\Phi$ is the *gauge function* defined by

$$\Phi(\lambda) := \frac{1}{(2\pi)^d}\int_{\mathbf{R}^d} \mathrm{e}^{-\lambda\Psi(\xi)}\,\mathrm{d}\xi \qquad \text{for all } \lambda > 0. \tag{1.5}$$

A3. For all $a > 0$, there exists a constant $A_a > 0$ such that

$$\Psi(a\xi) \ge A_a\Psi(\xi) \qquad \text{for all } \xi \in \mathbf{R}^d. \tag{1.6}$$

***Remark 1.1.*** In order to better understand condition A3, we note that $\Phi$ is non-increasing. Therefore, in particular, $\Phi(2\lambda) \le \Phi(\lambda)$ for all $\lambda > 0$. Condition A3 implies that a converse holds. Namely,

$$\limsup_{\lambda \downarrow 0} \frac{\Phi(\lambda)}{\Phi(2\lambda)} < \infty. \tag{1.7}$$

Unfortunately, (1.7) by itself does not seem to be enough to imply our main result.

Given an analytic Euclidean set $A$, we let $\dim_{\mathrm{H}} A$ denote its Hausdorff dimension [21, Appendix C], with the proviso added that the statement "$\dim_{\mathrm{H}} A < 0$" means that $A$ is empty.

We are ready to present the main result of this paper; the remainder of the article is dedicated to proving this fact.

**Theorem 1.2.** *Under conditions* A1–A3*, the following are equivalent:*

(1) *Almost surely, $u(t,x) = 0$ for some $t > 0$ and $x \in \mathbf{R}$;*
(2) *with positive probability, $u(t,x) = 0$ for some $t > 0$ and $x \in \mathbf{R}$;*
(3) $\int_0^1 \lambda \Phi(\lambda)\, \mathrm{d}\lambda < \infty$.

*In addition, if one of these conditions holds, then*

$$\dim_{\mathrm{H}} u^{-1}\{0\} = 2 - \overline{\mathrm{ind}}\Phi \qquad \text{a.s., where } \overline{\mathrm{ind}}\Phi := \limsup_{\lambda \downarrow 0} \frac{\log \Phi(\lambda)}{\log(1/\lambda)}. \tag{1.8}$$

In order to better understand Theorem 1.2, let us emphasize the special case where $\Psi(\xi) = \frac{1}{2}\|\xi\|^\alpha$ for all $\xi \in \mathbf{R}^d$, and $\alpha \in (0,2]$ is a fixed 'index of stability'. In this case, $\dot{L}$ is called the *isotropic stable Lévy noise* with index $\alpha$ and it is easy to see that $\Phi(\lambda) = \mathrm{const} \cdot \lambda^{d/\alpha}$. Hence, Theorem 1.2 yields the following result.

**Corollary 1.3.** *Consider the solution $u$ to* (1.1)*, where $\dot{L}$ is the isotropic stable Lévy noise with index $\alpha \in (0,2]$. Then, $u$ has zeros if and only if $d < 2\alpha$. Moreover, if $d < 2\alpha$, then $\dim_{\mathrm{H}} u^{-1}\{0\} = 2 - (d/\alpha)$ a.s.*

When $\alpha = 2$, $\dot{L}$ is the standard $d$-dimensional white noise on $\mathbf{R}^2$ and (1.1) simplifies to the more common form of the linear stochastic wave equation with $d$ independent, identically distributed components. In that case, there is indeed a large literature on this topic; see, for example, Dalang [8, 9], Dalang and Walsh [11], Cabaña [3, 4], Carmona and Nualart [6, 7], Gaveau [15, 16], Pyasetskaya [32] and Walsh [35].

One could also consider the case where $\Psi(\xi) = \xi' \boldsymbol{Q} \xi$, $\boldsymbol{Q} := (\boldsymbol{Q}_{ij})$ being a non-singular $d \times d$ covariance matrix. This leads to the weakly interacting system

$$\frac{\partial^2 u_i}{\partial t^2}(t,x) = \frac{\partial^2 u_i}{\partial x^2}(t,x) + \sum_{j=1}^d \boldsymbol{Q}_{ij} \dot{W}_j(t,x), \qquad t \geq 0, x \in \mathbf{R}, 1 \leq i \leq d, \tag{1.9}$$

with the initial condition $u_i(x,0) = \partial_t u_i(x,0) = 0$. Here, $\dot{W} := (\dot{W}_1, \ldots, \dot{W}_d)$ denotes a standard $d$-dimensional white noise on $\mathbf{R}^2$. Corollary 1.3 holds in this case with $\alpha = 2$, the result being that the solution has zeros if and only if $d < 4$. This result is closely related to a hyperbolic problem that was solved in Dalang and Nualart [10]. If we apply their characterization to the problem of the existence of zeros, then we obtain the same critical dimension of four as we do here. We add that the theorem of Dalang and Nualart [10] also holds in the presence of smooth multiplicative nonlinearities.

The present paper is organized as follows. In Section 2 we recall Lévy sheets and their associated Lévy processes, and present a brief introduction to stochastic integration with respect to a general Lévy noise. Sections 3, 4 and 5 are devoted to the study of the polar sets of the zero set of a Lévy sheet and of the random field defined as a sum of a Lévy sheet and an additive Lévy process, in a general setting. We finally prove Theorem 1.2 in Section 6.

In the rest of this section, we give some notation that will be used throughout. The $k$-dimensional Lebesgue measure on $\mathbf{R}^k$ is denoted by $\lambda_k$. For all $x \in \mathbf{R}^k$, $\|x\| := (x_1^2 + \cdots + x_k^2)^{1/2}$ denotes the $\ell^2$-norm of $x$ and $|x| := |x_1| + \cdots + |x_k|$ its $\ell^1$-norm.

The underlying parameter space is $\mathbf{R}^N$, or $\mathbf{R}_+^N = [0, \infty)^N$. A typical parameter $t \in \mathbf{R}^N$ is written as $t = (t_1, \ldots, t_N)$. There is a natural partial order $\preceq$ on $\mathbf{R}^N$. Namely, $s \preceq t$ if and only if $s_i \leq t_i$ for all $i = 1, \ldots, N$. When it is the case that $s \preceq t$, we define the interval $[s, t] = \prod_{i=1}^N [s_i, t_i]$. Finally, $s \wedge t$ denotes the $N$-vector whose $i$th coordinate is the minimum of $s_i$ and $t_i$.

We denote by $\hat{f}$ the normalized Fourier transform of $f \in L^1(\mathbf{R}^d)$, given by

$$\hat{f}(\xi) = \int_{\mathbf{R}^d} \mathrm{e}^{\mathrm{i}x \cdot \xi} f(x)\, \mathrm{d}x \qquad \text{for all } \xi \in \mathbf{R}^d. \tag{1.10}$$

## 2. Lévy sheets and their associated processes

In this section, we study the structure of the distribution function of the noise $\dot{L}$, which is called a Lévy sheet. In fact, we proceed by first studying such sheets under greater generality than is needed for (1.1). This is primarily because we view Lévy sheets as fundamental objects themselves.

### 2.1. Lévy sheets

Let us begin by introducing some notation for $N$-parameter, $d$-dimensional Lévy sheets. For more detailed information, particularly in the two-parameter case, consult Dalang and Walsh [11], Section 2.

Let $\dot{L}$ be a totally scattered symmetric $d$-dimensional random measure on $\mathbf{R}^N$ that is infinitely divisible in the following sense:

1. if $A$ and $B$ are disjoint Borel subsets of $\mathbf{R}^N$, then $\dot{L}(A)$ and $\dot{L}(B)$ are independent;
2. for every Borel set $A \subset \mathbf{R}^N$ of finite Lebesgue measure, the law of $\dot{L}(A)$ is described by

   $$\mathrm{E}[\exp\{\mathrm{i}\xi \cdot \dot{L}(A)\}] = \exp(-\lambda_N(A)\Psi(\xi)) \qquad \text{for all } \xi \in \mathbf{R}^d, \tag{2.1}$$

   where $\Psi$ is a real (and hence non-negative) negative definite function in the sense of Schoenberg [34]; see also the following remark.

***Remark 2.1.*** We recall Schoenberg's theorem: $\Psi$ is *negative definite* if and only if $\Psi(0) \geq 0$ and $\xi \mapsto \exp(-t\Psi(\xi))$ is positive definite in the sense of Herzog and Bochner [1, Theorem 7.8, page 41]. Equivalently, $\Psi$ is negative definite if it satisfies the Lévy–Khintchine formula [2, 33]. In the present case, $\Psi$ is also real-valued and, therefore, $\Psi(\xi) \geq 0$ for all $\xi \in \mathbf{R}^d$.

The *Lévy sheet* $L$ with Lévy exponent $\Psi$ is the $N$-parameter $d$-dimensional random field $\{L(t); t \in \mathbf{R}_+^N\}$ defined as

$$L(t) := \dot{L}([0,t]) \qquad \text{for all } t := (t_1, \ldots, t_N) \in \mathbf{R}_+^N. \tag{2.2}$$

Next, we mention a few commonly used families of Lévy noises and their corresponding Lévy sheets.

***Example 2.2 (White and stable noises).*** Choose and fix a constant $\chi \in (0, \infty)$. When $\Psi(\xi) = \chi \|\xi\|^2$ for all $\xi \in \mathbf{R}^d$, $\dot{L}$ is called the $d$-dimensional *white noise* on $\mathbf{R}^N$ and $L$ is called the $N$-parameter, $d$-dimensional *Brownian sheet*. The noise $\dot{L}$ is the (usual) *standard white noise* when $\chi = 1/2$ and the random field $L$ is then called the (usual) *standard Brownian sheet*.

More generally, if $\Psi(\xi) = \chi \|\xi\|^\alpha$ for $\xi \in \mathbf{R}^d$, then, by the Lévy–Khintchine formula, $\alpha \in (0, 2]$ and $\dot{L}$ is called the $d$-dimensional isotropic *stable noise* on $\mathbf{R}^N$ with index $\alpha$. In this case, $L$ is called the $N$-parameter *isotropic stable sheet* in $\mathbf{R}^d$ with index $\alpha$; see Ehm [14].

There are other interesting stable noises that are symmetric. For instance, one could consider $\Psi(\xi) = \frac{1}{2} \sum_{j=1}^d |\xi_j|^\alpha$. This gives rise to the symmetric noise each of whose $d$ components are independent, identically distributed one-dimensional stable noises on $\mathbf{R}^N$ with common index $\alpha$. The resulting random field is a stable sheet with independent, identically distributed coordinates, each with the same stability index $\alpha$.

***Example 2.3 (Noises with stable components).*** Let $\alpha_1, \ldots, \alpha_d \in (0, 2]$ be fixed, and consider $\Psi(\xi) = \chi \sum_{j=1}^d |\xi_j|^{\alpha_j}$, where $\chi \in (0, \infty)$ is fixed. The resulting random noise $\dot{L}$ has *stable components* with *index* $(\alpha_1, \ldots, \alpha_d)$ and the corresponding random field $L$ is the *stable sheet with stable components*. These were introduced in the one-parameter case by Pruitt [31] and investigated further by Hendricks [17, 18, 19, 20].

Let us also mention the following, immediate consequence of Theorem 1.2.

**Corollary 2.4.** *Consider the solution $u$ to (1.1), where $\dot{L}$ is Lévy noise with stable components with index $(\alpha_1, \ldots, \alpha_d) \in (0, 2]^d$. Then, $u$ has zeros if and only if $\sum_{j=1}^d (1/\alpha_j) < 2$. If and when $\sum_{j=1}^d (1/\alpha_j) < 2$, then*

$$\dim_{\mathrm{H}} u^{-1}\{0\} = 2 - \sum_{j=1}^d \frac{1}{\alpha_j} \qquad a.s. \tag{2.3}$$

## 2.2. The associated processes

Because $\Psi$ is a Lévy exponent, there exists a $d$-dimensional Lévy process $X := \{X(t); t \geq 0\}$ whose law is described uniquely by $\mathrm{E}[\mathrm{e}^{\mathrm{i}\xi\cdot X(t)}] = \mathrm{e}^{-t\Psi(\xi)}$ for all $t \geq 0$ and $\xi \in \mathbf{R}^d$. We may refer to $X$ as the *Lévy process associated to* $\dot{L}$.

By the inversion theorem, the transition densities of $X$ are given by

$$f(t;x) = \frac{1}{(2\pi)^d}\int_{\mathbf{R}^d} \cos(x\cdot\xi)\mathrm{e}^{-t\Psi(\xi)}\mathrm{d}\xi \qquad \text{for all } t \geq 0, x \in \mathbf{R}^d. \tag{2.4}$$

In particular, $\Phi(\lambda)$ is none other than $f(\lambda;0)$.

Now, let $X_1, \ldots, X_N$ be $N$ i.i.d. copies of $X$. Then, the *additive Lévy process associated to* $\dot{L}$ is defined as the $d$-dimensional $N$-parameter random field $\mathfrak{X} := \{\mathfrak{X}(t); t \in \mathbf{R}^N_+\}$, where

$$\mathfrak{X}(t) := X_1(t_1) + \cdots + X_N(t_N) \qquad \text{for all } t \in \mathbf{R}^N_+. \tag{2.5}$$

The density function of $\mathfrak{X}(t)$ at $0 \in \mathbf{R}^d$ is $\Phi(\sum_{i=1}^N t_i)$. This should be contrasted with the fact that the density function of $L(t)$ at $0 \in \mathbf{R}^d$ is $\Phi(\prod_{i=1}^N t_i)$.

## 2.3. Stochastic integrals

In this section, we proceed to construct Wiener-type stochastic integrals of the type $\dot{L}(\varphi) := \int_{\mathbf{R}^N} \varphi(t)\dot{L}(\mathrm{d}t)$, where $\varphi : \mathbf{R}^N \to \mathbf{R}$ is non-random, measurable, bounded and compactly supported. There exist integration theories that construct $\dot{L}(\varphi)$ abstractly; see, for example, the Bartle-type integrals of Dunford and Schwartz [13]. But, in order to mention some of the probabilistic properties of these integrals, we opt for a more direct approach that is closer to the original method of Wiener. As our method uses standard ideas to produce such stochastic integrals, we will only sketch the main steps.

First, suppose $\varphi$ is a simple function. That is, $\varphi(t) = \sum_{j=1}^m c_j \mathbf{1}_{A_j}$, where $A_1, \ldots, A_m$ are disjoint Borel sets in $\mathbf{R}^N$ and $c_1, \ldots, c_m$ are real constants. For such $\varphi$, we define $\dot{L}(\varphi) := \sum_{j=1}^m c_j \dot{L}(A_j)$. Because $\Psi(0) = 0$,

$$\begin{aligned}\mathrm{E}[\exp\{\mathrm{i}\xi\cdot\dot{L}(\varphi)\}] &= \exp\left(-\sum_{j=1}^m \lambda^N(A_j)\Psi(c_j\xi)\right)\\ &= \exp\left(-\int_{\mathbf{R}^N} \Psi(\varphi(t)\xi)\,\mathrm{d}t\right).\end{aligned} \tag{2.6}$$

Next, we consider a bounded, compactly-supported measurable function $\varphi : \mathbf{R}^N \to \mathbf{R}$. Standard measure-theoretic arguments reveal that we can find simple functions $\varphi_1, \varphi_2, \ldots : \mathbf{R}^N \to \mathbf{R}_+$ and a compact set $K \subset \mathbf{R}^N$ such that: (i) $\lim_{n\to\infty} \varphi_n = \varphi$ pointwise; (ii) $R := \sup_{n\geq 1}\sup_{t\in\mathbf{R}^N} |\varphi_n(t)| < \infty$; and (iii) $\varphi_n(t) = 0$ for all $t \notin K$ and all $n \geq 1$.

By linearity, and since $\Psi$ is real and nonnegative,

$$\begin{aligned} 1 &\geq \liminf_{n,m\to\infty} \mathrm{E}[\exp(\mathrm{i}\xi\cdot\{\dot{L}(\varphi_n)-\dot{L}(\varphi_m)\})] \\ &\geq \exp\left(-\limsup_{n,m\to\infty}\int_K \Psi(\{\varphi_n(t)-\varphi_m(t)\}\xi)\,\mathrm{d}t\right). \end{aligned} \tag{2.7}$$

Let $c := R\max_{1\leq j\leq d}|\xi_j|$ and note that $\Psi$ is uniformly continuous on $[-c,c]^d$. Therefore, $\lim_{n,m\to\infty}\mathrm{E}[\exp(\mathrm{i}\xi\cdot\{\dot{L}(\varphi_n)-\dot{L}(\varphi_m)\})]=1$. Consequently, the sequence $\{\dot{L}(\varphi_n)\}_{n\geq 1}$ is Cauchy in $L^0(\mathrm{P})$, and hence $\dot{L}(\varphi_n)$ converges in probability. The limit is denoted by $\dot{L}(\varphi)$ and has the following properties. As they can be proven directly from the preceding construction, we list them below without proof.

**Proposition 2.5.** *Let $\mathscr{B}$ denote the algebra of all non-random, bounded, compactly supported measurable functions from $\mathbf{R}^N$ to $\mathbf{R}$. Then, there is an 'iso-Lévy' process $\{\dot{L}(\varphi)\}_{\varphi\in\mathscr{B}}$ with the following properties:*

(i) $\mathrm{E}[\exp\{\mathrm{i}\xi\cdot\dot{L}(\varphi)\}]=\exp(-\int_{\mathbf{R}^N}\Psi(\varphi(t)\xi)\,\mathrm{d}t)$ *for all $\xi\in\mathbf{R}^d$ and $\varphi\in\mathscr{B}$;*
(ii) *if $\varphi_1,\varphi_2\in\mathscr{B}$, then $\dot{L}(\varphi_1+\varphi_2)=\dot{L}(\varphi_1)+\dot{L}(\varphi_2)$ a.s.;*
(iii) *if $a\in\mathbf{R}$ and $\varphi\in\mathscr{B}$, then $\dot{L}(a\varphi)=a\dot{L}(\varphi)$ a.s.*

That is, $\dot{L}$ defines a *random linear functional*, in the sense of Minlos [26, 27]. Alternatively, one might write $\int\varphi\,\mathrm{d}L$ in place of $\dot{L}(\varphi)$.

## 3. The image of a Lévy sheet

Theorem 1.2 is a consequence of the more general results of this section. Before we describe them, we need to introduce some notation.

Throughout the rest of the paper, $\{L(t); t\in\mathbf{R}^N_+\}$ denotes a $d$-dimensional $N$-parameter Lévy sheet whose Lévy exponent $\Psi$ satisfies the symmetry condition A1 and whose gauge function $\Phi$ satisfies the regularity conditions A2 and A3. Occasionally, we need the corresponding Lévy noise, which we denote by $\dot{L}$. Note that the Lévy sheet that arises from the SPDE (1.1) is one such random field with $N=2$.

We assume that the underlying sample space $\Omega$ is the collection of all ($N$-parameter) cadlag functions $\omega:\mathbf{R}^N_+\mapsto\mathbf{R}^d$. That is, $\omega\in\Omega$ if and only if:

1. for any net of elements $\{t_\alpha\}_{\alpha\in\mathbf{R}}$ in $\mathbf{R}^N_+$ such that $t_\alpha\preceq t_\beta$ whenever $\alpha\leq\beta$, and $t_\alpha$ converges to some $t\in\mathbf{R}^N_+$ as $\alpha\to\infty$, we have $\lim_{\alpha\to\infty}\omega(t_\alpha)=\omega(t)$;
2. for any net of elements $\{t_\alpha\}_{\alpha\in\mathbf{R}}$ in $\mathbf{R}^N_+$ such that $t_\alpha\preceq t_\beta$ whenever $\alpha\geq\beta$, it is the case that $\lim_{\alpha\to\infty}\omega(t_\alpha)$ exists.

We say that the Lévy process $L$ is in *canonical form* if

$$L(t)(\omega)=\omega(t)\qquad\text{for all }\omega\in\Omega\text{ and }t\in\mathbf{R}^N_+. \tag{3.1}$$

Throughout, we will assume that our Lévy process $L$ is in canonical form under a fixed probability measure P. This assumption is made tacitly and does not incur any loss of generality.

Define $\mathrm{P}_x$ to be the law of the process $x+L$ for every $x\in\mathbf{R}^d$. That is, for every Borel subset $A$ of $\Omega$,

$$\mathrm{P}_x\{\omega\in\Omega:\omega\in A\}=\mathrm{P}\{\omega\in\Omega:x+\omega\in A\}. \tag{3.2}$$

Let $\mathrm{E}_x$ denote the corresponding expectation operator. We will be primarily interested in the $\sigma$-finite (but infinite) measure

$$\mathrm{P}_{\lambda_d}(\bullet):=\int_{\mathbf{R}^d}\mathrm{P}_x(\bullet)\,\mathrm{d}x. \tag{3.3}$$

We will write $\mathrm{E}_{\lambda_d}$ for the corresponding 'expectation operator'. That is,

$$\mathrm{E}_{\lambda_d}(Z):=\int_\Omega Z(\omega)\mathrm{P}_{\lambda_d}(\mathrm{d}\omega)=\int_{\mathbf{R}^d}\mathrm{E}_x(Z)\,\mathrm{d}x \qquad \text{for all } Z\in L^1(\mathrm{P}_{\lambda_d}). \tag{3.4}$$

Let $\mathscr{P}(G)$ denote the collection of all probability measures on any Euclidean set $G$ and define, for all $\mu\in\mathscr{P}(\mathbf{R}^N_+)$,

$$\mathrm{I}(\mu):=\iint\Phi(\|s-t\|)\mu(\mathrm{d}t)\mu(\mathrm{d}s). \tag{3.5}$$

For all compact sets $G\subset\mathbf{R}^N_+$, we also consider the *capacity* of $G$,

$$\mathrm{Cap}(G):=\left[\inf_{\mu\in\mathscr{P}(G)}\mathrm{I}(\mu)\right]^{-1}, \tag{3.6}$$

where $\inf\varnothing:=\infty$ and $1/\infty:=0$. Let $L(G):=\bigcup_{t\in G}\{L(t)\}$ denote the range of $G$ under the random map $t\mapsto L(t)$. The following is the main result of this section.

**Theorem 3.1.** *Let $G$ be a non-random compact subset of $(0,\infty)^N$. There then exists a positive and finite constant $c=c(G)$ such that*

$$c^{-1}\,\mathrm{Cap}(G)\le\mathrm{P}_{\lambda_d}\{0\in L(G)\}=\mathrm{E}[\lambda_d(L(G))]\le c\,\mathrm{Cap}(G). \tag{3.7}$$

*Moreover, $c$ depends on $G$ only through $\inf_{x\in G}|x|$ and $\sup_{x\in G}|x|$.*

As a consequence of this theorem and Theorem 5.1 of Khoshnevisan and Xiao [23], we have the following equivalence theorem between the Lévy sheets of this paper and their associated additive Lévy processes.

**Corollary 3.2.** *Let $\{\mathfrak{X}(t);t\in\mathbf{R}^N_+\}$ denote the additive Lévy process associated to $L$. Then, for all non-random compact sets $G$ in $\mathbf{R}^N_+$, there exists a positive and finite constant $c=c(G)$ such that*

$$c^{-1}\mathrm{E}[\lambda_d(\mathfrak{X}(G))]\le\mathrm{E}[\lambda_d(L(G))]\le c\mathrm{E}[\lambda_d(\mathfrak{X}(G))]. \tag{3.8}$$

*Moreover, c depends only on G through* $\inf_{x\in G}|x|$ *and* $\sup_{x\in G}|x|$.

***Example 3.3.*** Suppose $L$ is an isotropic $\alpha$-stable Lévy sheet. Then, we can either calculate directly, or apply the preceding together with Corollary 5.4 of Khoshnevisan and Xiao [23], to find that

$$\mathrm{E}[\lambda_d(L(G))] > 0 \quad \Leftrightarrow \quad \mathrm{Cap}(G) > 0 \quad \Leftrightarrow \quad \mathscr{C}_{d/\alpha}(G) > 0, \tag{3.9}$$

where $\mathscr{C}_\beta$ denotes the standard $\beta$-dimensional Riesz capacity; see [21], Appendix D. □

***Example 3.4.*** Suppose $L$ is a Lévy sheet with stable components of index $(\alpha_1, \dots, \alpha_d) \in [0,2)^d$. Then,

$$\mathrm{E}[\lambda_d(L(G))] > 0 \quad \Leftrightarrow \quad \mathrm{Cap}(G) > 0 \quad \Leftrightarrow \quad \mathscr{C}_{\sum_{j=1}^d (1/\alpha_j)}(G) > 0, \tag{3.10}$$

where $\mathscr{C}_\beta$ is as in the previous example.

Our proof of Theorem 3.1 proceeds by first establishing an elementary lemma.

**Lemma 3.5.** *For all* $f \in L^1(\mathbf{R}^d)$ *and* $t \in \mathbf{R}_+^N$,

$$\mathrm{E}_{\lambda_d}[f(L(t))] = \int_{\mathbf{R}^d} f(x)\,\mathrm{d}x. \tag{3.11}$$

**Proof.** Because $f$ is measurable and integrable, Fubini's theorem applies and we have

$$\mathrm{E}_{\lambda_d}[f(L(t))] = \mathrm{E}\left[\int_{\mathbf{R}^d} f(x + L(t))\,\mathrm{d}x\right], \tag{3.12}$$

this being manifestly equal to $\int_{\mathbf{R}^d} f(x)\,\mathrm{d}x$. □

The next lemma is nearly as simple, but particularly useful to us.

**Lemma 3.6.** *If f and g are two probability densities on* $\mathbf{R}^d$ *such that* $\hat{f}, \hat{g} \in L^1(\mathbf{R}^d)$, *then*

$$\mathrm{E}_{\lambda_d}[f(L(t))g(L(s))] = \frac{1}{(2\pi)^d}\int_{\mathbf{R}^d} \mathrm{e}^{-\ell(s,t)\Psi(\xi)}\hat{g}(\xi)\overline{\hat{f}(\xi)}\,\mathrm{d}\xi, \tag{3.13}$$

*where* $\ell$ *is the symmetric function*

$$\ell(s,t) := \lambda_N([0,s] \triangle [0,t]) \qquad \textit{for all } s,t \in \mathbf{R}_+^N. \tag{3.14}$$

**Proof.** We write

$$\mathrm{E}_{\lambda_d}[f(L(t))g(L(s))] = \mathrm{E}\left[\int_{\mathbf{R}^d} f(x + L(t))g(x + L(s))\,\mathrm{d}x\right]$$

$$= \int_{\mathbf{R}^d} f(y) \mathrm{E}[g(y + L(s) - L(t))]\, \mathrm{d}y \qquad (3.15)$$
$$= \int_{\mathbf{R}^d} f(y) \int_{\mathbf{R}^d} g(y-z) \mu_{s,t}(\mathrm{d}z)\, \mathrm{d}y,$$

where $\mu_{s,t}$ denotes the distribution of $L(t) - L(s)$. By the inversion theorem of Fourier analysis,

$$\begin{aligned} \mathrm{E}_{\lambda_d}[f(L(s))g(L(t))] &= \frac{1}{(2\pi)^d} \int_{\mathbf{R}^d} \mathrm{d}y\, f(y) \int_{\mathbf{R}^d} \mu_{s,t}(\mathrm{d}z) \int_{\mathbf{R}^d} \mathrm{e}^{-\mathrm{i}(y-z)\cdot\xi} \hat{g}(\xi)\, \mathrm{d}\xi \\ &= \frac{1}{(2\pi)^d} \int_{\mathbf{R}^d} \hat{\mu}_{s,t}(\xi) \hat{g}(\xi) \overline{\hat{f}(\xi)}\, \mathrm{d}\xi. \end{aligned} \qquad (3.16)$$

In order to compute the Fourier transform of $\mu_{s,t}$, we first note that

$$\begin{aligned} L(t) - L(s) &= \dot{L}([0,t]) - \dot{L}([0,s]) \\ &= \dot{L}([0,t] \setminus [0,s]) - \dot{L}([0,s] \setminus [0,t]). \end{aligned} \qquad (3.17)$$

The last two random variables are independent of one another. Moreover, by symmetry, $-\dot{L}([0,s] \setminus [0,t])$ has the same distribution as $\dot{L}([0,s] \setminus [0,t])$. Therefore, $L(t) - L(s)$ has the same distribution as $\dot{L}([0,s] \triangle [0,t])$. From this and (2.1), it follows that $\hat{\mu}_{s,t}(\xi) = \exp(-\ell(s,t)\Psi(\xi))$. We plug this into (3.16) in order to conclude the proof. □

Next, we recall a well-known estimate for the function $\ell$ that was defined in (3.14); see [21], Lemma 1.3.1, page 460, for a proof.

**Lemma 3.7.** *If $u \preceq v$ are in $(0,\infty)^N$, then there exist positive and finite constants $a$ and $b$ such that*

$$a\|s-t\| \le \ell(s,t) \le b\|s-t\| \qquad \text{for all } s,t \in [u,v]. \qquad (3.18)$$

**Lemma 3.8.** *If $G \subset \mathbf{R}_+^N$ is compact and non-random, then*

$$\mathrm{P}_{\lambda_d}\{0 \in \overline{L(G)}\} \ge \mathrm{const} \cdot \mathrm{Cap}(G). \qquad (3.19)$$

*Moreover, the constant depends on $G$ only through $\inf_{x\in G}|x|$ and $\sup_{x \in G}|x|$.*

**Proof.** Without loss of generality, we assume that $\mathrm{Cap}(G) > 0$; otherwise there is nothing to prove. In that case, there exists $\mu \in \mathscr{P}(G)$ such that

$$\mathrm{I}(\mu) \le \frac{2}{\mathrm{Cap}(G)} < \infty. \qquad (3.20)$$

We choose and fix this $\mu$ throughout.

For all probability density functions $f$ on $\mathbf{R}^d$, define the random variable

$$J(f;\mu) := \int_{\mathbf{R}_+^N} f(L(t))\mu(\mathrm{d}t). \tag{3.21}$$

By Lemma 3.5, $\mathrm{E}_{\lambda_d}[J(f;\mu)] = 1$. On the other hand, by Lemma 3.6,

$$\begin{aligned} \mathrm{E}_{\lambda_d}(|J(f;\mu)|^2) &= \frac{1}{(2\pi)^d} \iint_{(\mathbf{R}_+^N)^2} \mu(\mathrm{d}s)\mu(\mathrm{d}t) \int_{\mathbf{R}^d} \mathrm{e}^{-\ell(s,t)\Psi(\xi)} |\hat{f}(\xi)|^2 \, \mathrm{d}\xi \\ &\le \frac{1}{(2\pi)^d} \iint_{(\mathbf{R}_+^N)^2} \mu(\mathrm{d}s)\mu(\mathrm{d}t) \int_{\mathbf{R}^d} \mathrm{e}^{-\ell(s,t)\Psi(\xi)} \, \mathrm{d}\xi. \end{aligned} \tag{3.22}$$

We now apply Lemma 3.7 to see that there exists a positive and finite constant $a$ – depending on $G$ only through $\inf_{x\in G}|x|$ and $\sup_{x\in G}|x|$ – such that

$$\begin{aligned} \mathrm{E}_{\lambda_d}(|J(f;\mu)|^2) &\le \frac{1}{(2\pi)^d} \iint \Phi(a\|s-t\|)\mu(\mathrm{d}s)\mu(\mathrm{d}t) \\ &\le A_a \mathrm{I}(\mu) \le \frac{2A_a}{\mathrm{Cap}(G)}. \end{aligned} \tag{3.23}$$

The second bound follows from A3 and the third from (3.20).

By the Paley–Zygmund inequality, [21], Lemma 1.4.1, Chapter 3, can be applied to the $\sigma$-finite measure $\mathrm{P}_{\lambda_d}$ and this implies that

$$\mathrm{P}_{\lambda_d}\{J(f;\mu) > 0\} \ge \frac{|(\mathrm{E}_{\lambda_d}[J(f;\mu)]|^2}{\mathrm{E}_{\lambda_d}(|J(f;\mu)|^2)} \ge \frac{\mathrm{Cap}(G)}{2A_a}. \tag{3.24}$$

Let $f$ be a probability density function on $B(0;1)$, where $B(x;r)$ denotes the open (Euclidean) ball of radius $r > 0$ centered at $x \in \mathbf{R}^d$. For all $\delta > 0$, define $f_\delta(x) := \delta^{-d} f(x/\delta)$ as $x$ varies over $\mathbf{R}^d$. Note that $f_\delta$ is a probability density on $B(0;\delta)$. Furthermore, if $\omega$ is such that $J(f_\delta;\mu)(\omega) > 0$, then there exists $t \in G$ such that $L(t)(\omega) \in B(0;\delta)$. Hence, it follows that

$$\mathrm{P}_{\lambda_d}\{L(G) \cap B(0;\delta) \ne \varnothing\} \ge \frac{\mathrm{Cap}(G)}{2A_a}. \tag{3.25}$$

From here, we can deduce the lemma by letting $\delta$ tend to zero. □

For all $t \in \mathbf{R}_+^N$, we denote by $\mathscr{F}^0(t)$ the $\sigma$-algebra generated by $\{L(s); s \preceq t\}$. Let $\mathscr{F}^1(t)$ denote the completion of $\mathscr{F}^0(t)$ with respect to $\mathrm{P}_x$ for all $x \in \mathbf{R}^d$. Finally, define

$$\mathscr{F}(t) := \bigcap_{u \succeq t} \mathscr{F}^1(u). \tag{3.26}$$

In the language of the general theory of random fields, $\mathscr{F} := \{\mathscr{F}(t); t \in \mathbf{R}_+^N\}$ is the *augmented* $N$-parameter filtration of $L$.

The following is an analogue of the 'Markov property' in the present setting.

**Proposition 3.9.** *If $s \preceq t$ are in $\mathbf{R}_+^N$, then for all measurable functions $f : \mathbf{R}^d \mapsto \mathbf{R}_+$,*

$$\mathrm{E}_{\lambda_d}[f(L(t))|\mathscr{F}(s)] = (T_{s,t}^L f)(L(s)), \tag{3.27}$$

*where*

$$(T_{s,t}^L f)(y) := \mathrm{E}[f(L(t) - L(s) + y)] \qquad \text{for all } y \in \mathbf{R}_+^N. \tag{3.28}$$

**Proof.** Consider measurable functions $f, g, h_1, \ldots, h_m : \mathbf{R}^d \mapsto \mathbf{R}_+$ and times $t, s, s^1, \ldots, s^m \in \mathbf{R}_+^N$ such that $t \succeq s \succeq s^j$, for all $j = 1, \ldots, m$. Because $\dot{L}(A)$ and $\dot{L}(B)$ are independent when $A \cap B = \varnothing$,

$$\begin{aligned} &\mathrm{E}_{\lambda_d}\left[f(L(t)) \cdot g(L(s)) \cdot \prod_{j=1}^m h_j(L(s^j))\right] \\ &= \int_{\mathbf{R}^d} \mathrm{E}\left[f(L(t)+x) \cdot g(L(s)+x) \cdot \prod_{j=1}^m h_j(L(s^j)+x)\right] \mathrm{d}x \qquad (3.29) \\ &= \int_{\mathbf{R}^d} \mathrm{E}[f(L(t)-L(s)+y)]\mathrm{E}\left[\prod_{j=1}^m h_j(L(s^j)-L(s)+y)\right] g(y)\,\mathrm{d}y. \end{aligned}$$

Set $f = h_1 = \cdots = h_m \equiv 1$ to see that $\mathrm{P}_{\lambda_d}\{L(s) \in \bullet\} = \lambda_d(\bullet)$, and the desired result follows. □

The next result is Cairoli's maximal inequality. This inequality is proved as a consequence of the *commuting* property of the $N$-parameter filtration $\mathscr{F}$.

**Proposition 3.10.** *For all $Y \in L^2(\mathrm{P}_{\lambda_d})$, there exists a modification of $t \mapsto \mathrm{E}_{\lambda_d}[Y|\mathscr{F}(t)]$ that is right-continuous in every variable, uniformly in the remaining $N-1$ variables, and satisfies*

$$\mathrm{E}_{\lambda_d}\left(\sup_{t \in \mathbf{R}_+^N} |\mathrm{E}_{\lambda_d}[Y|\mathscr{F}(t)]|^2\right) \le 4^N \mathrm{E}_{\lambda_d}[Y^2]. \tag{3.30}$$

***Remark 3.11.*** It might help to recall that a function $F : \mathbf{R}_+^N \to \mathbf{R}$ is right-continuous in each variable, uniformly in all its other variables, if and only if for all $j = 1, \ldots, N$ and $t_j \in \mathbf{R}_+$,

$$\lim_{s_j \downarrow t_j} \sup |F(s) - F(t)| = 0, \tag{3.31}$$

where the supremum is taken over all $s_i = t_i \in \mathbf{R}_+$ $(1 \le i \ne j \le N)$.

**Proof of Proposition 3.10.** When $L$ is the Brownian sheet and the infinite measure $\mathrm{P}_{\lambda_d}$ is replaced by the probability measure P, this result follows from Cairoli's maximal inequality and the Cairoli–Walsh commutation theorem; see Chapter 7 of Khoshnevisan [21], Section 2.

In order to prove the present form of the proposition, we need to replace the Gaussian Brownian sheet with the more general random field $L$ and, more significantly, P with $\mathrm{P}_{\lambda_d}$. Fortunately, none of this requires us to introduce too many new ideas. Therefore, we merely outline the requisite changes in the existing theory to accomodate the present formulation of this proposition.

We appeal to a close analogy with the usual $N$-parameter theory of processes and say that $\mathscr{F}$ is a *commuting $N$-parameter filtration* with respect to $\mathrm{P}_{\lambda_d}$ if for all $s, t \in \mathbf{R}_+^N$ and every $\mathscr{F}(s)$-measurable random variable $Y \in L^2(\mathrm{P}_{\lambda_d})$,

$$\mathrm{E}_{\lambda_d}[Y | \mathscr{F}(t)] = \mathrm{E}_{\lambda_d}[Y | \mathscr{F}(s \wedge t)], \tag{3.32}$$

off a $\mathrm{P}_{\lambda_d}$-null set. In the case that $N = 2$ and $\mathrm{P}_{\lambda_d}$ is replaced by P, (3.32) essentially reduces to hypothesis (F4) of Cairoli and Walsh [5]. The general case [P] is studied – using the same methods as Cairoli and Walsh – in Chapter 7 of Khoshnevisan [21], Section 2.1, among other places. The same methods also prove (3.32) for $\mathrm{P}_{\lambda_d}$. We will not reproduce these steps, as they are standard.

It is known that (3.32) implies the 'Cairoli maximal inequality' announced in the statement of the proposition. In the case that $\mathrm{P}_{\lambda_d}$ is replaced everywhere by P, this is covered by Theorem 2.3.2 of [21], page 235. One adapts the proof so that the result continues to holds for $\mathrm{P}_{\lambda_d}$. This proves the proposition. □

**Lemma 3.12.** *If $G$ is a non-random compact set in $\mathbf{R}_+^N$, then*

$$\mathrm{P}_{\lambda_d}\{0 \in L(G)\} \leq \mathrm{const} \cdot \mathrm{Cap}(G). \tag{3.33}$$

*Moreover, the constant depends on $G$ only through the quantities $\sup_{x \in G} |x|$ and $\inf_{x \in G} |x|$.*

**Proof.** Without loss of generality, we may assume that

$$\mathrm{P}_{\lambda_d}\{0 \in L(G)\} > 0. \tag{3.34}$$

Otherwise, there is nothing to prove.

Let $f$ be a probability density on $\mathbf{R}^d$ and consider an arbitrary $\mu \in \mathscr{P}(G)$. In accordance with Proposition 3.9, for all $s \in G$,

$$\mathrm{E}_{\lambda_d}[J(f;\mu) | \mathscr{F}(s)] \geq \int_{t \succeq s} (T_{s,t}^L f)(L(s)) \mu(\mathrm{d}t) \qquad \mathrm{P}_{\lambda_d}\text{-a.s.}, \tag{3.35}$$

where $J(f;\mu)$ is defined in (3.21). The two sides of (3.35) are right-continuous in $s$, both coordinatewise; see Dalang and Walsh [11], and also Ehm [14]. Therefore, for all probability densities $f$ on $\mathbf{R}^d$ such that $\hat{f} \in L^1(\mathbf{R}^d)$, and for all $\mu \in \mathscr{P}(G)$, the following statement holds.

$$\text{There exists one null set off which (3.35) holds for all } s \in G. \tag{3.36}$$

Path regularity of $L$ implies the existence of a random variable $\sigma$ with values in $G \cup \{\rho\}$ – where $\rho$ is an abstract 'cemetery' point not in $G$ – which has the following properties:

1. $\sigma = \rho$ if and only if $L(s) \neq 0$ for every $s \in G$;
2. $L(\sigma) = 0$ on $\{\sigma \neq \rho\}$.

For all integers $k \geq 1$ and all Borel sets $E \subset \mathbf{R}^N_+$, define

$$\mu_k(E) := \frac{\mathrm{P}_{\lambda_d}\{\sigma \in E, \sigma \neq \rho, |L(0)| \leq k\}}{\mathrm{P}_{\lambda_d}\{\sigma \neq \rho, |L(0)| \leq k\}}. \tag{3.37}$$

Thanks to (3.34) and the monotone convergence theorem,

$$\lim_{k\to\infty} \mathrm{P}_{\lambda_d}\{\sigma \neq \rho, |L(0)| \leq k\} = \mathrm{P}_{\lambda_d}\{0 \in L(G)\} > 0. \tag{3.38}$$

Therefore, there exists $k_0 > 0$ such that $\mu_k$ is well defined for all $k \geq k_0$. In fact, $\mu_k$ is a probability measure on $G$, provided that $k \geq k_0$.

Define, for all $k \geq k_0$ and all probability densities $f$ on $\mathbf{R}^d$,

$$Q_k(f) := \sup_{u \in G} \mathrm{E}_{\lambda_d}[J(f;\mu_k)|\mathscr{F}(u)]. \tag{3.39}$$

Then, thanks to (3.35) and (3.36), and because $L(\sigma) = 0$ a.s. on $\{\sigma \neq \rho\}$,

$$\begin{aligned} Q_k(f) &\geq \mathrm{E}_{\lambda_d}[J(f;\mu_k)|\mathscr{F}(s)]|_{s=\sigma} \cdot \mathbf{1}_{\{\sigma\neq\rho\}} \\ &\geq \mathbf{1}_{\{\sigma\neq\rho, |L(0)|\leq k\}} \cdot \int_{t\succeq\sigma} (T^L_{\sigma,t} f)(0)\mu_k(\mathrm{d}t) \qquad \mathrm{P}_{\lambda_d}\text{-a.s.}, \end{aligned} \tag{3.40}$$

provided that $\hat{f} \in L^1(\mathbf{R}^d)$. (Specifically, (3.36) ensures that we can apply (3.35) with the choice $s := \sigma(\omega)$ for every $\omega$ such that $\sigma(\omega) \neq \rho$.) We can square both sides and integrate $[\mathrm{P}_{\lambda_d}]$ to find that

$$\begin{aligned} &\mathrm{E}_{\lambda_d}(|Q_k(f)|^2) \\ &\geq \mathrm{P}_{\lambda_d}\{\sigma \neq \rho, |L(0)| \leq k\} \cdot \int_{\mathbf{R}^N_+} \left( \int_{t\succeq s} (T^L_{s,t} f)(0)\mu_k(\mathrm{d}t) \right)^2 \mu_k(\mathrm{d}s). \end{aligned} \tag{3.41}$$

Since $\mu_k$ is a probability measure on $G$ $[k \geq k_0]$, the Cauchy–Schwarz inequality implies that

$$\mathrm{E}_{\lambda_d}(|Q_k(f)|^2) \geq \Lambda^2 \mathrm{P}_{\lambda_d}\{\sigma \neq \rho, |L(0)| \leq k\}, \tag{3.42}$$

where

$$\Lambda := \int_{\mathbf{R}_+^N} \int_{t \succeq s} (T_{s,t}^L f)(0) \mu_k(\mathrm{d}t) \mu_k(\mathrm{d}s). \tag{3.43}$$

Suppose also that $\hat{f}(\xi) \geq 0$ for all $\xi \in \mathbf{R}^d$. Because $f, \hat{f} \in L^1(\mathbf{R}^d)$, we can apply the inversion theorem to $f$ – in the same manner as in the proof of Lemma 3.6 – to find that for all $s, t \in G$ with $s \preceq t$,

$$\begin{aligned} (T_{s,t}^L f)(0) &= \frac{1}{(2\pi)^d} \int_{\mathbf{R}^d} \mathrm{e}^{-\ell(s,t)\Psi(\xi)} \hat{f}(\xi) \, \mathrm{d}\xi \\ &\geq \frac{1}{(2\pi)^d} \int_{\mathbf{R}^d} \mathrm{e}^{-b\|s-t\|\Psi(\xi)} \hat{f}(\xi) \, \mathrm{d}\xi. \end{aligned} \tag{3.44}$$

Here, $b$ is a positive and finite constant that depends only on the distance between $G$ and the axes of $\mathbf{R}_+^N$, as well as $\sup_{x \in G} |x|$; see Lemma 3.7. Consequently,

$$\Lambda \geq \frac{1}{2^N (2\pi)^d} \iint_{(\mathbf{R}_+^N)^2} \mu_k(\mathrm{d}s) \mu_k(\mathrm{d}t) \int_{\mathbf{R}^d} \mathrm{e}^{-b\|s-t\|\Psi(\xi)} \hat{f}(\xi) \, \mathrm{d}\xi. \tag{3.45}$$

According to A3, $b\Psi(\xi) \leq \Psi(\xi / A_{1/b})$. Therefore, the preceding implies that

$$\Lambda \geq \frac{A_{1/b}^d}{2^N (2\pi)^d} \iint_{(\mathbf{R}_+^N)^2} \mu_k(\mathrm{d}s) \mu_k(\mathrm{d}t) \int_{\mathbf{R}^d} \mathrm{e}^{-\|s-t\|\Psi(\xi)} \hat{f}(\xi A_{1/b}) \, \mathrm{d}\xi. \tag{3.46}$$

Now, thanks to Proposition 3.10 and (3.22),

$$\begin{aligned} \mathrm{E}_{\lambda_d}(|Q_k(f)|^2) &\leq 4^N \mathrm{E}_{\lambda_d}(|J(f;\mu_k)|^2) \\ &\leq \frac{4^N}{(2\pi)^d} \iint_{(\mathbf{R}_+^N)^2} \mu_k(\mathrm{d}s) \mu_k(\mathrm{d}t) \int_{\mathbf{R}^d} \mathrm{e}^{-\ell(s,t)\Psi(\xi)} |\hat{f}(\xi)|^2 \, \mathrm{d}\xi. \end{aligned}$$

Therefore, Lemma 3.7 implies that there exists a positive and finite constant $a$ that depends only on the distance between $G$ and the axes of $\mathbf{R}_+^N$, as well as $\sup_{x \in G} |x|$, such that

$$\begin{aligned} \mathrm{E}_{\lambda_d}(|Q_k(f)|^2) &\leq \frac{4^N}{(2\pi)^d} \iint_{(\mathbf{R}_+^N)^2} \mu_k(\mathrm{d}s) \mu_k(\mathrm{d}t) \int_{\mathbf{R}^d} \mathrm{e}^{-a\|s-t\|\Psi(\xi)} |\hat{f}(\xi)|^2 \, \mathrm{d}\xi \\ &\leq \frac{4^N}{(2\pi)^d A_a^d} \iint_{(\mathbf{R}_+^N)^2} \mu_k(\mathrm{d}s) \mu_k(\mathrm{d}t) \int_{\mathbf{R}^d} \mathrm{e}^{-\|s-t\|\Psi(\xi)} |\hat{f}(\xi / A_a)|^2 \, \mathrm{d}\xi; \end{aligned}$$

see A3 for the last inequality. This, (3.42) and (3.46) together imply that

$$c \iint_{(\mathbf{R}_+^N)^2} \mu_k(\mathrm{d}s) \mu_k(\mathrm{d}t) \int_{\mathbf{R}^d} \mathrm{e}^{-\|s-t\|\Psi(\xi)} |\hat{f}(\xi / A_a)|^2 \, \mathrm{d}\xi$$

$$\geq \left( \iint_{(\mathbf{R}_+^N)^2} \mu_k(\mathrm{d}s)\mu_k(\mathrm{d}t) \int_{\mathbf{R}^d} \mathrm{e}^{-\|s-t\|\Psi(\xi)} \hat{f}(\xi A_{1/b})\, \mathrm{d}\xi \right)^2 \qquad (3.47)$$

$$\times \mathrm{P}_{\lambda_d}\{\sigma \neq \rho, |L(0)| \leq k\},$$

where $c := 16^N (2\pi)^d / (A_a A_{1/b}^2)^d$.

The preceding is valid for every probability density $f$ on $\mathbf{R}^d$ whose Fourier transform is integrable and nonnegative. We now make a particular choice of $f$. Namely, we replace $f$ by $f_\eta$, where $f_\eta(x) := \eta^{-d} g(x/\eta)$, for $\eta > 0$, and $g$ is the density function of an isotropic stable random variable of index $\alpha \in (0, 2]$ satisfying

$$\frac{1}{\alpha} \geq \log_2(A_a A_{1/b}). \qquad (3.48)$$

Note that $\hat{g}(\xi) = \exp(-\|\xi\|^\alpha)$ for all $\xi \in \mathbf{R}^d$. Hence, as $\hat{f}_\eta(\xi) = \exp(-\eta^\alpha \|\xi\|^\alpha)$ for all $\xi \in \mathbf{R}^d$, we have

$$|\hat{f}_\eta(\xi/A_a)|^2 = \exp\left( -\frac{2\eta^\alpha \|\xi\|^\alpha}{A_a^\alpha} \right) \leq \exp(-\eta^\alpha \|\xi\|^\alpha A_{1/b}^\alpha) = \hat{f}_\eta(\xi A_{1/b}). \qquad (3.49)$$

Thus, we find that

$$\begin{aligned} &\mathrm{P}_{\lambda_d}\{\sigma \neq \rho, |L(0)| \leq k\} \\ &\leq c \left[ \iint_{\|s-t\| \geq \theta} \mu_k(\mathrm{d}s)\mu_k(\mathrm{d}t) \int_{\mathbf{R}^d} \mathrm{e}^{-\|s-t\|\Psi(\xi)} \exp\left( -\frac{2\eta^\alpha \|\xi\|^\alpha}{A_a^\alpha} \right) \mathrm{d}\xi \right]^{-1}, \end{aligned} \qquad (3.50)$$

where $\theta > 0$ is an arbitrary parameter. Consequently, for all $k \geq k_0$,

$$\begin{aligned} &\mathrm{P}_{\lambda_d}\{0 \in L(G), |L(0)| \leq k\} \\ &\leq \tilde{c} \left[ \iint_{\|s-t\| \geq \theta} \Phi(\|s-t\|) \mu_k(\mathrm{d}s)\mu_k(\mathrm{d}t) \right]^{-1}. \end{aligned} \qquad (3.51)$$

We may observe that $\{\mu_k\}_{k \geq k_0}$ is a collection of probability measures, all of which are supported on the same compact set $G$. Therefore, by Prohorov's theorem, we can extract a subsequence of $k$'s along which $\mu_k$ converges weakly to some $\mu_* \in \mathscr{P}(G)$. Because $\Phi$ is uniformly continuous on $[\theta, \infty)$, it follows that

$$\begin{aligned} &\liminf_{k \to \infty} \mathrm{P}_{\lambda_d}\{0 \in L(G), |L(0)| \leq k\} \\ &\leq \tilde{c} \left[ \iint_{\|s-t\| \geq \theta} \Phi(\|s-t\|) \mu_*(\mathrm{d}s)\mu_*(\mathrm{d}t) \right]^{-1}. \end{aligned} \qquad (3.52)$$

By the monotone convergence theorem, the left-hand side is precisely $\mathrm{P}_{\lambda_d}\{0 \in L(G)\}$. Finally, let $\theta \downarrow 0$ and appeal to the monotone convergence theorem once again to find that $\mathrm{P}_{\lambda_d}\{0 \in L(G)\} \leq \tilde{c}/\mathrm{I}(\mu_*)$. This is at most $\tilde{c}\,\mathrm{Cap}(G)$ by the definition of $\mathrm{Cap}(G)$. This concludes the proof. □

We conclude this section by proving Theorem 3.1.

**Proof of Theorem 3.1.** First note that

$$\begin{aligned} \mathrm{P}_{\lambda_d}\{0 \in L(G)\} &= \int_{\mathbf{R}^d} \mathrm{P}\{-x \in L(G)\}\,\mathrm{d}x \\ &= \mathrm{E}[\lambda_d(L(G))]. \end{aligned} \tag{3.53}$$

This proves the identity in Theorem 3.1. Moreover, Lemma 3.12 proves that the preceding is at most $c\,\mathrm{Cap}(G)$, whence follows the upper bound of the theorem.

Similarly,

$$\mathrm{P}_{\lambda_d}\{0 \in \overline{L(G)}\} = \mathrm{E}[\lambda_d(\overline{L(G)})]. \tag{3.54}$$

But the set-difference between $L(G)$ and its Euclidean closure is a denumerable union of sets, each of which is at most $(d-1)$-dimensional. See Dalang and Walsh [11], Proposition 2.1, for the case $N=2$; the general case is proved similarly. It follows that

$$\mathrm{E}[\lambda_d(\overline{L(G)} \setminus L(G))] = 0. \tag{3.55}$$

Thus,

$$\begin{aligned} \mathrm{P}_{\lambda_d}\{0 \in L(G)\} &= \mathrm{E}[\lambda_d(L(G))] \\ &= \mathrm{E}[\lambda_d(\overline{L(G)})] \\ &= \mathrm{P}_{\lambda_d}\{0 \in \overline{L(G)}\}. \end{aligned} \tag{3.56}$$

Apply Lemma 3.8 to complete the proof. □

# 4. A perturbed random field

Let $\mathfrak{X}$ be a $d$-dimensional $N$-parameter additive Lévy process with Lévy exponent $\nu\Psi$, where $\nu \in (0,\infty)$ is fixed. We also assume that $\mathfrak{X}$ is independent of $L$. Define

$$F(t) := L(t) + \mathfrak{X}(t) \qquad \text{for all } t \in \mathbf{R}_+^N. \tag{4.1}$$

The following proves that the contents of Theorem 3.1 remain unchanged if $L$ is replaced by the perturbed random field $F$.

**Theorem 4.1.** *Let $G$ denote a non-random compact subset of $(0,\infty)^N$. There then exists a positive and finite constant $c = c(G,\nu)$ such that*

$$c^{-1}\,\mathrm{Cap}(G) \leq \mathrm{P}_{\lambda_d}\{0 \in F(G)\} = \mathrm{E}[\lambda_d(F(G))] \leq c\,\mathrm{Cap}(G). \tag{4.2}$$

*Moreover, $c$ depends on $G$ only through* $\inf_{x\in G}|x|$ *and* $\sup_{x\in G}|x|$.

**Proof.** Write $\mathrm{P}_z$ for the law of $z+F$, and $\mathrm{P}_{\lambda_d} := \int_{\mathbf{R}^d} \mathrm{P}_z\,\mathrm{d}z$, as before. Also, $\mathrm{E}_{\lambda_d}$ denotes the corresponding expectation operator. As in Lemma 3.5, if $f \in L^1(\mathbf{R}^d)$ and $t \in \mathbf{R}_+^N$, then

$$\mathrm{E}_{\lambda_d}[f(F(t))] = \int_{\mathbf{R}^d} f(x)\,\mathrm{d}x. \tag{4.3}$$

Also, the proof of Lemma 3.6 implies that if $f$ and $g$ are two probability densities on $\mathbf{R}^d$ such that $\hat{f}, \hat{g} \in L^1(\mathbf{R}^d)$, then for all $s, t \in \mathbf{R}_+^N$,

$$\mathrm{E}_{\lambda_d}[f(F(t))g(F(s))] = \frac{1}{(2\pi)^d}\int_{\mathbf{R}^d} \mathrm{e}^{-(\ell(s,t)+c\sum_{j=1}^N |t_j-s_j|)\Psi(\xi)}\,\hat{g}(\xi)\overline{\hat{f}(\xi)}\,\mathrm{d}\xi. \tag{4.4}$$

Finally, as in Proposition 3.9, we have

$$\mathrm{E}_{\lambda_d}[f(F(t))|\mathscr{F}(s)] = (T_{s,t}^F f)(F(s)), \tag{4.5}$$

valid for all $s, t \in \mathbf{R}_+^N$ and all measurable functions $f : \mathbf{R}^d \to \mathbf{R}_+$. Here,

$$(T_{s,t}^F f)(y) := \mathrm{E}[f(F(t)+y)] \qquad \text{for all } s,t \in \mathbf{R}_+^N, y \in \mathbf{R}^d, \tag{4.6}$$

and $\mathscr{F}$ denotes now the filtration of $F$. As the filtration of $\mathfrak{X}$ is also commuting [23], Lemma 4.2 and its proof, the remainder of the proof of Theorem 3.1 goes through in the present case without any drastic changes. □

## 5. Polar sets for the level sets

The aim of this section is to compute the Hausdorff dimension of the zero set $L^{-1}\{0\}$ of an $N$-parameter, $d$-dimensional Lévy sheet $L$, provided that Conditions A1, A2 and A3 are assumed throughout.

Recall that a non-random set $G$ is *polar* for a random set $S$ if $\mathrm{P}\{S \cap G \neq \varnothing\} = 0$. We begin by describing all sets that are polar for the level sets of the random field $L$ and its perturbation $F$.

**Theorem 5.1.** *Let $G$ be a non-random compact subset of $(0,\infty)^N$. The following are then equivalent:*

(1) $\mathrm{P}\{L^{-1}\{0\} \cap G \neq \varnothing\} > 0$;
(2) $\mathrm{P}\{F^{-1}\{0\} \cap G \neq \varnothing\} > 0$;
(3) $\mathrm{Cap}(G) > 0$.

We can apply the preceding, together with Corollary 2.13 of Khoshnevisan and Xiao [23], in order to deduce the following equivalence theorem between the Lévy sheets of this paper and their associated additive Lévy processes.

**Corollary 5.2.** *Let $\{\mathfrak{X}(t)\}_{t\in\mathbf{R}_+^N}$ denote the additive Lévy process associated to $L$. Then, for all non-random compact sets $G$ in $(0,\infty)^N$, there exists a positive and finite constant $c$ such that*

$$c^{-1}\mathrm{P}\{\mathfrak{X}^{-1}\{0\}\cap G\neq\varnothing\}\leq \mathrm{P}\{L^{-1}\{0\}\cap G\neq\varnothing\}\leq c\mathrm{P}\{\mathfrak{X}^{-1}\{0\}\cap G\neq\varnothing\}.$$

*Moreover, $c$ depends on $G$ only through $\inf_{x\in G}|x|$ and $\sup_{x\in G}|x|$.*

**Proof of Theorem 5.1.** There exists $\alpha>0$ such that $G\subseteq(\alpha,1/\alpha)^N$. Therefore, if $t\in G$, then $\mathrm{P}_x$-a.s. for all $x\in\mathbf{R}^d$,

$$L(t)=L(a)+\dot{L}([0,t]\setminus[0,a]), \tag{5.1}$$

where $a:=(\alpha,\ldots,\alpha)\in(0,\infty)^N$ and the two terms on the right-hand side describe independent random fields. We can decompose the rightmost term as

$$\dot{L}([0,t]\setminus[0,a])=\sum_{j=1}^N L^j(\pi^j(t-a))+\tilde{L}(t-a), \tag{5.2}$$

where:

(1) $L^1,\ldots,L^N,\tilde{L}$ are totally independent $d$-dimensional random fields;
(2) each $L^j$ is an $(N-1)$-parameter Lévy sheet with Lévy exponent $\Psi$;
(3) $\tilde{L}$ is an $N$-parameter Lévy sheet with Lévy exponent $\Psi$;
(4) $\pi^j$ maps $t\in\mathbf{R}_+^N$ to $\pi^j t\in\mathbf{R}_+^{N-1}$, which is the same vector as $t$, but with $t_j$ removed.

This identity is valid $\mathrm{P}_x$-a.s. for all $x\in\mathbf{R}^d$. But, because we can choose path-regular versions of $L^1,\ldots,L^N$ and $\tilde{L}$, it follows that $\mathrm{P}_x$-a.s. for all $x\in\mathbf{R}^d$,

$$L(t)=L(a)+\sum_{j=1}^N L^j(\pi^j(t-a))+\tilde{L}(t-a)\qquad\text{for all } t\in G. \tag{5.3}$$

Moreover, the $N+2$ processes on the right-hand side (viewed as random fields indexed by $t$) are totally independent under $\mathrm{P}_x$ for all $x\in\mathbf{R}^d$. Note that

$$\begin{aligned}
&\mathrm{P}\{L^{-1}\{0\}\cap G\neq\varnothing\}\\
&\quad=\mathrm{P}\{\exists t\in G: L(t)=0\}\\
&\quad=\int_{\mathbf{R}^d}\mathrm{P}\left\{\exists t\in G:\sum_{j=1}^N L^j(\pi^j(t-a))+\tilde{L}(t-a)=-x\right\}\mathrm{P}\{L(a)\in\mathrm{d}x\}\\
&\quad=\int_{\mathbf{R}^d}\mathrm{P}\left\{\exists t\in G:\sum_{j=1}^N L^j(\pi^j(t-a))+L(t-a)=-x\right\}\mathrm{P}\{L(a)\in\mathrm{d}x\}.
\end{aligned}\tag{5.4}$$

Thanks to (1.5) and the inversion theorem,

$$\frac{\mathrm{P}\{L(a)\in \mathrm{d}x\}}{\mathrm{d}x}=\frac{1}{(2\pi)^d}\int_{\mathbf{R}^d}\exp(-\mathrm{i}x\cdot\xi-\alpha^N\Psi(\xi))\,\mathrm{d}\xi \tag{5.5}$$

and this is bounded above, uniformly for all $x\in\mathbf{R}^d$, by the positive and finite constant $c:=c(\alpha):=(2\pi)^{-d}\int_{\mathbf{R}^d}\exp(-\alpha^N\Psi(\xi))\,\mathrm{d}\xi$. Consequently,

$$\begin{aligned}&\mathrm{P}\{L^{-1}\{0\}\cap G\neq\varnothing\}\\ &\quad\le c\int_{\mathbf{R}^d}\mathrm{P}\left\{\exists t\in G:\sum_{j=1}^N L^j(\pi^j(t-a))+L(t-a)=-x\right\}\mathrm{d}x\\ &\quad= c\mathrm{E}[\lambda_d(W(G\ominus\{a\}))],\end{aligned} \tag{5.6}$$

where $W:=\sum_{j=1}^N(L^j\circ\pi^j)+L$ and $G\ominus\{a\}:=\{t-a;t\in G\}$. Because $G\ominus\{a\}$ is a non-random compact subset of $(0,\infty)^N$, one can prove – in exactly the same manner that we proved Theorem 4.1 – that $\mathrm{E}[\lambda_d(W(G\ominus\{a\}))]\le\mathrm{const}\cdot\mathrm{Cap}(G\ominus\{a\})$, where the constant depends only on $G$ and $\alpha$. We omit the details. Because $\mathrm{I}(\mu)$ is convolution-based, $G\mapsto\mathrm{Cap}(G)$ is translation invariant, whence we have $\mathrm{Cap}(G\ominus\{a\})=\mathrm{Cap}(G)$. This proves that (1) implies (3).

Next, we prove that (3) implies (1). As above, let $\alpha>0$ be such that $G\subset(\alpha,1/\alpha)^N$. For any $\epsilon>0$ and $\mu\in\mathscr{P}(\mathbf{R}_+^N)$, we define a random measure on $\mathbf{R}_+^N$ by

$$\mathscr{J}_\epsilon(B)=\frac{1}{(2\epsilon)^d}\int_B \mathbf{1}_{B(0,\epsilon)}(L(t))\mu(\mathrm{d}t)\qquad\text{for all Borel sets } B\subset\mathbf{R}_+^N. \tag{5.7}$$

We shall need the following two lemmas.

**Lemma 5.3.** *For every $\mu\in\mathscr{P}(G)$, $\liminf_{\epsilon\to0^+}\mathrm{E}[\mathscr{J}_\epsilon(G)]>0$.*

**Proof.** By Fatou's lemma,

$$\begin{aligned}\liminf_{\epsilon\to0^+}\mathrm{E}[\mathscr{J}_\epsilon(G)]&=\liminf_{\epsilon\to0^+}\int_G\frac{\mathrm{P}\{|L(t)|\le\epsilon\}}{(2\epsilon)^d}\mu(\mathrm{d}t)\\ &\ge\int_G p(t)\mu(\mathrm{d}t),\end{aligned} \tag{5.8}$$

where $p(t)$ denotes the density function of $L(t)$ at zero. By the inversion theorem, $p(t)=\Phi(\prod_{j=1}^N t_j)$ and this is strictly positive uniformly for all $t\in G$. This implies the desired result. □

**Lemma 5.4.** *Let $K:\mathbf{R}_+^N\times\mathbf{R}_+^N\mapsto\mathbf{R}_+$ be a measurable function. There then exists a finite positive constant $c$, depending only on $G$, such that for every $\mu\in\mathscr{P}(G)$ and all*

$\epsilon > 0$,

$$\mathrm{E}\left[\iint_{G\times G} K(s,t)\,\mathscr{J}_\epsilon(\mathrm{d}s)\,\mathscr{J}_\epsilon(\mathrm{d}t)\right] \leq c \iint_{G\times G} K(s,t)\Phi(\|s-t\|)\mu(\mathrm{d}t)\mu(\mathrm{d}s).$$

*In particular,*

$$\mathrm{E}(|\mathscr{J}_\epsilon(G)|^2) \leq c\mathrm{I}(\mu). \tag{5.9}$$

**Proof.** By Fubini's theorem,

$$\begin{aligned}&\mathrm{E}\left[\iint_{G\times G} K(s,t)\,\mathscr{J}_\epsilon(\mathrm{d}s)\,\mathscr{J}_\epsilon(\mathrm{d}t)\right]\\ &\quad = \frac{1}{(2\epsilon)^{2d}} \iint_{G\times G} K(s,t)\mathrm{P}\{|L(s)|\leq\epsilon, |L(t)|\leq\epsilon\}\mu(\mathrm{d}t)\mu(\mathrm{d}s).\end{aligned} \tag{5.10}$$

We define

$$\begin{aligned} N_1(s,t) &:= L(s) - L(s\wedge t),\\ N_2(s,t) &:= L(t) - L(s\wedge t).\end{aligned} \tag{5.11}$$

Clearly,

$$\mathrm{P}\{|L(s)|\leq\epsilon, |L(t)|\leq\epsilon\} \leq \mathrm{P}\{|L(s\wedge t) + N_1(s,t)|\leq\epsilon, |L(t)-L(s)|\leq 2\epsilon\}.$$

Note that $L(s\wedge t)$, $N_1(s,t)$ and $N_2(s,t)$ are mutually independent; this follows immediately by considering the representation of $L$ via $\dot{L}$. Thus,

$$\mathrm{P}\{|L(s)|\leq\epsilon, |L(t)|\leq\epsilon\} \leq \sup_{x\in\mathbf{R}^d} \mathrm{P}_x\{|L(s\wedge t)|\leq\epsilon\} \times \mathrm{P}\{|L(t)-L(s)|\leq 2\epsilon\}.$$

Because of symmetry, $L(s\wedge t)$ is a $16^d$-weakly unimodal random vector [23], Corollary 2.2. That is,

$$\sup_{x\in\mathbf{R}^d} \mathrm{P}_x\{|L(s\wedge t)|\leq\epsilon\} \leq 16^d \mathrm{P}\{|L(s\wedge t)|\leq\epsilon\}. \tag{5.12}$$

But, if $s,t\in G$, then $s\wedge t$ is also in $[\alpha, 1/\alpha]^N$. Hence, by (5.5), the density of $L(s\wedge t)$ is bounded, uniformly in $s$ and $t$ in $G$. Consequently, we can find a positive and finite constant $c_2 = c_2(\alpha)$ such that

$$\sup_{x\in\mathbf{R}^d} \mathrm{P}_x\{|L(s\wedge t)|\leq\epsilon\} \leq c_2(2\epsilon)^d \qquad \text{for all } \epsilon>0, s,t\in G. \tag{5.13}$$

Consequently,

$$\mathrm{E}\left[\iint_{G\times G} K(s,t)\,\mathscr{J}_\epsilon(\mathrm{d}s)\,\mathscr{J}_\epsilon(\mathrm{d}t)\right] \leq c_3 \iint_{G\times G} K(s,t)p(s,t)\mu(\mathrm{d}t)\mu(\mathrm{d}s),$$

where $p(s,t)$ is the density of $L(t)-L(s)$ at zero. By the inversion theorem,

$$\begin{aligned} p(s,t) &= \frac{1}{(2\pi)^d}\int_{\mathbf{R}^d} \mathrm{e}^{-\ell(s,t)\Psi(\xi)}\,\mathrm{d}\xi \\ &\le \frac{1}{(2\pi)^d}\int_{\mathbf{R}^d} \mathrm{e}^{-a\|s-t\|\Psi(\xi)}\,\mathrm{d}\xi, \end{aligned} \tag{5.14}$$

where $a\in(0,\infty)$ depends only on $G$; see Lemma 3.7. Thanks to A3, we can write

$$\begin{aligned} p(s,t) &\le \frac{1}{(2\pi)^d}\int_{\mathbf{R}^d} \mathrm{e}^{-\|s-t\|\Psi(A_a\xi)}\,\mathrm{d}\xi \\ &= \frac{1}{A_a^d}\Phi(\|t-s\|), \end{aligned} \tag{5.15}$$

and this concludes the proof. □

Let us now continue the proof that (3) implies (1). Note that (3) implies the existence of $\mu\in\mathscr{P}(G)$ such that $\mathrm{I}(\mu)<\infty$. There then exists a continuous function $\rho:\mathbf{R}^N\mapsto[1,\infty)$ such that $\lim_{s\mapsto s_0}\rho(s)=\infty$ for every $s_0\in\mathbf{R}^N$ with at least one coordinate equalling 0 and such that

$$\iint_{G\times G}\rho(s-t)\Phi(\|s-t\|)\mu(\mathrm{d}t)\mu(\mathrm{d}s)<\infty; \tag{5.16}$$

see Khoshnevisan and Xiao [22], page 73, for a construction of $\rho$.

Consider now the sequence of random measures $\{\mathscr{J}_\epsilon\}_{\epsilon>0}$. It follows from Lemmas 5.3, 5.4 and a second-moment argument that we can extract a subsequence $\{\mathscr{J}_{\epsilon_n}\}_{n\in\mathbf{N}}$ – converging weakly to a random measure $\mathscr{J}$ – such that

$$\mathrm{P}\{\mathscr{J}(G)>0\}\ge\frac{(\inf_{0<\epsilon<1}\mathrm{E}[\mathscr{J}_\epsilon(G)])^2}{\sup_{\epsilon>0}\mathrm{E}(|\mathscr{J}_\epsilon(G)|^2)}>0; \tag{5.17}$$

see Khoshnevisan and Xiao [22], proof of Lemma 3.6, and Khoshnevisan, Xiao and Shieh [25], page 26. Moreover,

$$\iint_{G\times G}\rho(s-t)\mathscr{J}(\mathrm{d}t)\mathscr{J}(\mathrm{d}s)\le c\iint_{G\times G}\rho(s-t)\Phi(\|s-t\|)\mu(\mathrm{d}t)\mu(\mathrm{d}s). \tag{5.18}$$

This and (5.16) together imply that a.s.,

$$\mathscr{J}\{s\in G: s_j=a \text{ for some } j\}=0 \qquad \text{for all } a\in\mathbf{R}_+. \tag{5.19}$$

In order to deduce (1), it suffices to prove that the random measure $\mathscr{J}$ is supported on $L^{-1}\{0\}\cup G$. For this purpose, it suffices to prove that for every $\delta>0$, $\mathscr{J}(D(\delta))=0$ a.s., where $D(\delta):=\{s\in G:|L(s)|>\delta\}$. Because of (5.19), the proof of this fact follows exactly the arguments of Khoshnevisan and Xiao [22], page 76.

The equivalence of (2) and (3) is proved in exactly the same way and is therefore omitted. □

Next, we derive our first dimension result.

**Theorem 5.5.** *Let $G$ be a non-random compact subset of $\mathbf{R}_+^N$. Then, almost surely on* $\{L^{-1}\{0\} \cap G \neq \varnothing\}$,

$$\dim_{\mathrm{H}}(L^{-1}\{0\} \cap G) = \sup\{q \in (0, N) : \mathrm{I}^{(q)}(\mu) < \infty \text{ for some } \mu \in \mathscr{P}(G)\},$$

*where*

$$\mathrm{I}^{(q)}(\mu) := \iint \frac{\Phi(\|t-s\|)}{\|t-s\|^q} \mu(\mathrm{d}s)\mu(\mathrm{d}t). \tag{5.20}$$

**Corollary 5.6.** *Almost surely on* $\{L^{-1}\{0\} \neq \varnothing\}$, $\dim_{\mathrm{H}} L^{-1}\{0\} = N - \overline{\mathrm{ind}}\Phi$.

***Example 5.7.*** If $L$ is an isotropic stable Lévy sheet of index $\alpha \in (0, 2]$, then we can combine Corollary 5.6 with Example 3.5 of Khoshnevisan, Xiao and Shieh [25] to find that a.s. on $\{L^{-1}\{0\} \cap G \neq \varnothing\}$,

$$\dim_{\mathrm{H}}(L^{-1}\{0\} \cap G) = \dim_{\mathrm{H}} G - \frac{d}{\alpha}. \tag{5.21}$$

Also, $\dim_{\mathrm{H}} L^{-1}\{0\} = N - (d/\alpha)$ a.s. on $\{L^{-1}\{0\} \neq \varnothing\}$. This last fact was first proven in [14].

It is not easy to produce a direct proof of Theorem 5.5. Therefore, we instead use a 'comparison principle' to deduce the theorem from an analogous (known) result on the associated additive Lévy process.

**Proof of Theorem 5.5.** Let $\mathfrak{X}$ denote the additive Lévy process associated to $L$ and choose and fix an integer $M \geq 1$ and a real number $\alpha \in (0, 2]$. Consider $M$ independent isotropic stable Lévy processes of index $\alpha$, $S^1, \ldots, S^M$, all taking values in $\mathbf{R}^N$. We assume that $\{S^1, \ldots, S^M, L, \mathfrak{X}\}$ are independent. Consider the $M$-parameter, $N$-valued additive stable Lévy process,

$$\mathfrak{S}(u) := S^1(u_1) + \cdots + S^M(u_M) \qquad \text{for all } u \in \mathbf{R}_+^M. \tag{5.22}$$

Suppose $N > \alpha M$. Then, by Propositions 4.7 and 4.8 of Khoshnevisan, Shieh and Xiao [25], we deduce that

$$\mathrm{P}\{\mathfrak{X}^{-1}\{0\} \cap G \cap \mathfrak{S}(\mathbf{R}_+^M) = \varnothing\} > 0 \quad \Leftrightarrow \quad \inf_{\mu \in \mathscr{P}(G)} \mathrm{I}^{(N-\alpha M)}(\mu) < \infty. \tag{5.23}$$

A standard approximation argument can be used in conjunction with Theorem 5.1 to deduce that for all non-random analytic sets $F \subset \mathbf{R}_+^N$,

$$\mathrm{P}\{\mathfrak{X}^{-1}\{0\} \cap F = \varnothing\} > 0 \quad \Leftrightarrow \quad \mathrm{P}\{L^{-1}\{0\} \cap F = \varnothing\} > 0. \tag{5.24}$$

Thanks to the independence properties of $\mathfrak{S}$, we can apply the preceding with $F := G \cap \mathfrak{S}(\mathbf{R}_+^M)$ (first condition on $\mathfrak{S}$ and then take expectations). Consequently, (5.23) implies that

$$\mathrm{P}\{L^{-1}\{0\} \cap G \cap \mathfrak{S}(\mathbf{R}_+^M) = \varnothing\} > 0 \quad \Leftrightarrow \quad \inf_{\mu \in \mathscr{P}(G)} \mathrm{I}^{(N-\alpha M)}(\mu) < \infty. \tag{5.25}$$

The rest proof of the proof of Theorem 5.5 now follows exactly as the proof of Theorem 3.2 of Khoshnevisan and Xiao [23] and is therefore omitted. □

We use a similar comparison principle to deduce Corollary 5.6.

**Proof of Corollary 5.6.** Let $\mathfrak{X}$ denote the additive Lévy process associated to $L$ and assume that $L$ and $\mathfrak{X}$ are independent. Thanks to Theorem 5.5, and Theorem 3.2 of Khoshnevisan and Xiao [25], the following holds with probability one:

$$\dim_{\mathrm{H}}(L^{-1}\{0\} \cap G)\mathbf{1}_{\{L^{-1}\{0\} \cap G \neq \varnothing\}} = \dim_{\mathrm{H}}(\mathfrak{X}^{-1}\{0\} \cap G)\mathbf{1}_{\{\mathfrak{X}^{-1}\{0\} \cap G \neq \varnothing\}}.$$

Now, consider a non-random *upright cube* $G \subset (0,\infty)^N$, that is, a closed set $G$ of the form $G = \prod_{i=1}^N [a_i, b_i]$. The proof of Theorem 1.1 of [25] then states that $\dim_{\mathrm{H}}(\mathfrak{X}^{-1}\{0\} \cap G) = N - \overline{\mathrm{ind}}\Phi$ almost surely on $\{\mathfrak{X}^{-1}\{0\} \cap G \neq \varnothing\}$.

Hence, $\dim_{\mathrm{H}}(L^{-1}\{0\} \cap G) = N - \overline{\mathrm{ind}}\Phi$ almost surely on $\{L^{-1}\{0\} \cap G \neq \varnothing\}$. Because we need to consider only upright cubes that have rational coordinates, a limiting argument finishes the proof. □

## 6. Proof of Theorem 1.2

We now return to the problem in the Introduction and study the SPDE (1.1).

Let $x_0 \in \mathbf{R}$ be fixed and let $\mathscr{S}(x_0)$ be the collection of all $(t,x) \in \mathbf{R}_+ \times \mathbf{R}$ such that $t \geq |x - x_0|$. Elementary geometric considerations lead us to the following equivalent formulation:

$$\mathscr{S}(x_0) = \{(t,x) \in \mathbf{R}_+ \times \mathbf{R} : (x_0, 0) \in \mathscr{C}(t,x)\}. \tag{6.1}$$

We will need the following well-known result; see, for example, the Introduction of Dalang and Walsh [11] for the case $x_0 = 0$.

**Lemma 6.1.** *Let $x_0 \in \mathbf{R}$ be fixed. We can then write*

$$u(t,x) = \tfrac{1}{2}\tilde{F}(t-x, t+x) \qquad \textit{for all } (t,x) \in \mathscr{S}(x_0), \tag{6.2}$$

*where $\tilde{F}$ is a copy of the perturbed random field $F$ of Section 4 with $N = 2$ and $\nu = \frac{1}{2}$.*

Our proof of Theorem 1.2 requires two more simple lemmas.

**Lemma 6.2.** *Let $x_0 \in \mathbf{R}$ be fixed and let $G$ be a non-random compact subset of $(0,\infty)^2 \cap \mathscr{S}(x_0)$. Then,*

$$\mathrm{P}\{u^{-1}\{0\} \cap G \neq \varnothing\} > 0 \quad \Leftrightarrow \quad \mathrm{Cap}(G) > 0. \tag{6.3}$$

**Proof.** Define the rotation map $\rho : \mathbf{R}^2 \to \mathbf{R}^2$ by $\rho(t,x) := (t-x, t+x)$. Its inverse is described by $\rho^{-1}(u,v) = \frac{1}{2}(v+u, v-u)$. According to Lemma 6.1,

$$u = \tfrac{1}{2}(\tilde{F} \circ \rho) \qquad \text{on } \mathscr{S}(x_0). \tag{6.4}$$

Now the process $\frac{1}{2}F$ is the perturbed version of $\frac{1}{2}L$ and the latter is a Lévy sheet with Lévy exponent $\xi \mapsto \Psi(\xi/2)$. Thanks to Theorem 5.1, the restriction of $F^{-1}\{0\}$ to $\mathscr{S}(x_0)$ has the same polar sets as the restriction of $L^{-1}\{0\}$ to $\mathscr{S}(x_0)$. Consequently, if $G$ is any non-random compact subset of $\mathscr{S}(x_0) \cap (0,\infty)^2$, then

$$\mathrm{P}\{u^{-1}\{0\} \cap G \neq \varnothing\} > 0 \quad \Leftrightarrow \quad \mathrm{Cap}(\rho^{-1}G) > 0. \tag{6.5}$$

Note that if $\mu \in \mathscr{P}(G)$, then $\mu \circ \rho \in \mathscr{P}(\rho^{-1}G)$ and $\mathrm{I}(\mu) = \mathrm{I}(\mu \circ \rho)$. Consequently, $\mathrm{Cap}(\rho^{-1}G) = \mathrm{Cap}(G)$ and the lemma follows. □

**Lemma 6.3.** *Let $G$ be a non-random compact subset of $(0,\infty)^2$. Then,*

$$\mathrm{P}\{u^{-1}\{0\} \cap G \neq \varnothing\} > 0 \quad \Leftrightarrow \quad \mathrm{Cap}(G) > 0. \tag{6.6}$$

**Proof.** If $\mathrm{P}\{u^{-1}\{0\} \cap G \neq \varnothing\}$ is positive, then there must exist $x_0 \in \mathbf{R}$ such that $\mathrm{P}\{u^{-1}\{0\} \cap G \cap \mathscr{S}(x_0) \neq \varnothing\}$ is positive. Lemma 6.2 then implies that $\mathrm{Cap}(G \cap \mathscr{S}(x_0)) > 0$. Because $\mathrm{Cap}(G) \geq \mathrm{Cap}(G \cap \mathscr{S}(x_0))$, this proves the '$\Rightarrow$' portion of the lemma.

The converse is proved similarly. Indeed, it can be shown that Cap is a Choquet capacity, whence $\mathrm{Cap}(G) = \sup_{x_0 \in \mathbf{R}} \mathrm{Cap}(G \cap \mathscr{S}(x_0))$, thanks to Choquet's capacitability theorem [12], Theorem 28, page 52-III. □

Finally, we prove Theorem 1.2 via a comparison argument.

**Proof of Theorem 1.2.** Let $\mathfrak{X}$ be the 2-parameter $d$-dimensional additive Lévy process that is associated to $L$. Lemma 6.3, used in conjunction with Theorem 5.1 and Corollary 5.2, implies that the zero sets of $u$ and $\mathfrak{X}$ have the same polar sets. Consequently, for all integers $n \geq 1$,

$$\begin{aligned} &\mathrm{P}\{u^{-1}\{0\} \cap [-n,n] \times [0,n] \neq \varnothing\} > 0 \\ &\quad \Leftrightarrow \quad \mathrm{P}\{\mathfrak{X}^{-1}\{0\} \cap [-n,n] \times [0,n] \neq \varnothing\} > 0. \end{aligned} \tag{6.7}$$

Let $n$ tend to infinity to find that

$$\mathrm{P}\{u^{-1}\{0\} \neq \varnothing\} > 0 \quad \Leftrightarrow \quad \mathrm{P}\{\mathfrak{X}^{-1}\{0\} \neq \varnothing\} > 0. \tag{6.8}$$

According to Proposition 3.1 of [22], this last condition holds if and only if

$$\int_{[0,1]^2} \Phi(\|z\|)\,\mathrm{d}z < \infty. \tag{6.9}$$

Integration in polar coordinates proves '(2) ⇔ (3)' in Theorem 1.2. Since '(1) ⇒ (2)' holds tautologically, in order to prove that (1) and (3) are equivalent, it suffices to prove that if $u^{-1}\{0\}$ is non-empty with positive probability, then it is non-empty a.s.

With this in mind, let us assume that $\mathrm{P}\{u^{-1}\{0\} \neq \varnothing\} > 0$. There then exists an integer $n \geq 1$ large enough such that

$$\mathrm{P}\{u^{-1}\{0\} \cap \mathscr{C}(n,0) \neq \varnothing\} > 0. \tag{6.10}$$

Thanks to (1.3), $u(t,x) = \frac{1}{2}\dot{L}(\mathscr{C}(t,x))$. We observed in Proposition 2.5 that $\dot{L}$ has the following two properties: (i) for all $a \in \mathbf{R}^2$ and $A \subset \mathbf{R}^2$ Borel-measurable, $\dot{L}(a+A)$ has the same law as $\dot{L}(A)$; (ii) $\dot{L}(A)$ and $\dot{L}(B)$ are independent whenever $A$ and $B$ are disjoint measurable subsets of $\mathbf{R}^2$. It follows from this that $\{u^{-1}\{0\} \cap \mathscr{C}(n, 2kn)\}_{k=0}^{\infty}$ are i.i.d. random sets. Used in conjunction with (6.10) and the Borel–Cantelli lemma for independent events, this implies that

$$\mathrm{P}\left\{ u^{-1}\{0\} \cap \bigcup_{k=0}^{\infty} \mathscr{C}(n, 2kn) \neq \varnothing \right\} = 1. \tag{6.11}$$

Because $\bigcup_{k=0}^{\infty} \mathscr{C}(n,2kn) \subset \mathbf{R}_+ \times \mathbf{R}$, this proves that (2) ⇒ (1).

In order to complete the proof of Theorem 1.2, we assume that $\int_0^1 \lambda\Phi(\lambda)\,\mathrm{d}\lambda$ is finite and proceed to compute the Hausdorff dimension of $u^{-1}\{0\}$. A comparison argument, similar to the one employed at the beginning of this proof, shows that if $G$ is a non-random compact subset of $(0,\infty)^2$, then a.s.,

$$\dim_{\mathrm{H}}(u^{-1}\{0\} \cap G)\mathbf{1}_{\{u^{-1}\{0\} \cap G \neq \varnothing\}} = \dim_{\mathrm{H}}(\mathfrak{X}^{-1}\{0\} \cap G)\mathbf{1}_{\{\mathfrak{X}^{-1}\{0\} \cap G \neq \varnothing\}}.$$

Let $G = [0,n] \times [-n,n]$ and then let $n \to \infty$ to find that a.s.,

$$\dim_{\mathrm{H}}(u^{-1}\{0\})\mathbf{1}_{\{u^{-1}\{0\} \neq \varnothing\}} = \dim_{\mathrm{H}}(\mathfrak{X}^{-1}\{0\})\mathbf{1}_{\{\mathfrak{X}^{-1}\{0\} \neq \varnothing\}}. \tag{6.12}$$

Because we have assumed that $\int_0^1 \lambda\Phi(\lambda)\,\mathrm{d}\lambda < \infty$, the already proven equivalence of (1)–(3) shows that $\mathbf{1}_{\{u^{-1}\{0\} \neq \varnothing\}} = 1$ a.s. Therefore, Theorem 1.1 of Khoshnevisan, Shieh and Xiao [25], (1.11), implies the stated formula for the Hausdorff dimension of $u^{-1}\{0\}$. □

## Acknowledgements

We thank an anonymous referee for making a number of helpful suggestions and for pointing out some errors and misprints. Research supported in part by the NSF Grant DMS-07-06728.